\documentclass[11pt,letterpaper]{article}

\usepackage{etex}
\usepackage{amsmath}
\usepackage{amsfonts}
\usepackage{amssymb}
\usepackage{multirow}
\usepackage{multicol}
\usepackage{longtable}
\usepackage[usenames,dvipsnames]{color}
\usepackage[top=1in, bottom=2in, left=1in, right=1in]{geometry} 
\usepackage{pdflscape}
\usepackage{graphicx}
\usepackage{setspace}
\usepackage{palatino}
\usepackage{url}
\usepackage[linesnumbered,algoruled]{algorithm2e}
\usepackage[pdftex,breaklinks]{hyperref}

\newcommand{\bit}{\begin{itemize}}
\newcommand{\eit}{\end{itemize}}
\newcommand{\ben}{\begin{enumerate}}
\newcommand{\een}{\end{enumerate}}

\newcommand{\T}{\rule[2.6ex]{0pt}{0pt}}
\newcommand{\B}{\rule[-1.2ex]{0pt}{0pt}}

\hypersetup{
    unicode=false,          
    pdftoolbar=true,        
    pdfmenubar=true,        
    pdffitwindow=false,     
    pdfstartview={Fit},     
	pdfpagelayout={SinglePage},		
    pdftitle={},    							
    pdfauthor={Anshul Agarwal},		
    pdfsubject={Documentation},		
    pdfcreator={Anshul Agarwal},	
    pdfproducer={Anshul Agarwal},	
    pdfkeywords={}, 							
    pdfnewwindow=true,      
    colorlinks=true,        
    linkcolor=blue,         
    citecolor=blue,         
    filecolor=blue,         
    urlcolor=blue           
}

\author{Anshul Agarwal\thanks{Corresponding author, anshul.agarwal@united.com, 
United Airlines, Chicago, IL 60606, USA}}
\date{}
\title{Multi-echelon Supply Chain Inventory Planning  using Simulation-Optimization with Data Resampling}

\begin{document}

\maketitle

\begin{abstract}

Modeling and optimization of multi-echelon supply chain systems is challenging as it requires a holistic approach that exploits synergies and interactions between echelons while accurately accounting for variability observed by these systems.  We develop a simulation-optimization framework that minimizes average inventory while maintaining desired average $\beta$ service level at stocking locations.  We use a discrete-event simulation framework to accurately capture system interactions.  Instead of a parametric estimation approach, the demand and the lead time variability are quantified by bootstrap sampling the historical data, thus preserving the true nature of the variability.  We compare three different open source simulation-optimization libraries - the derivative free methods from SciPy.Optimize, a Bayesian optimization algorithm Scikit-Optimize, and a radial basis function based black-box optimizer RBFOpt.  The experiments demonstrate practical applicability of our approach.  While we observe substantially lower inventory levels and computationally superior results from RBFOpt, depending on the problem, search strategy, and the random start states, close enough good solution can be obtained from both RBFOpt and Scikit-Optimize.  The optimization results demonstrate a preference for a centralized inventory planning scheme that help with risk pooling.  Moreover, with no order placement cost, the optimal solution tends to order more frequently in order to lower inventory.

\end{abstract}

\section{Introduction}

Complex end-to-end supply chains utilize multiple tiers or \emph{echelons} of stocking locations to not only minimize logistics costs, but also reduce inventory costs through risk pooling and efficient placement of inventory between echelons \cite{chopra2007supply, ganeshan1999managing, evant2003, zipkin2000foundations}.  New shipments are first stored at a central facility, which then acts as an internal supplier to other stocking locations as the product travels through multiple stages before reaching the final customer.  In this work we address the challenging problem of optimizing inventory in such complex systems.

Inventory optimization of a multi-echelon system can be approached by applying single-echelon optimization techniques to every echelon.  However, such an approach fails to achieve true network inventory optimization.  Disregarding the network view of inventory usage and applying stand-alone replenishment strategies to one echelon without considering impact on other echelons can lead to redundant safety stock or customer service failure despite adequate inventory \cite{evant2003}.

A multi-echelon setting can lead to sophisticated complex interactions between the inventory levels of upstream and downstream stocking locations.  Inventory unavailability at the upstream node can increase lead time for downstream facilities.  Centralized upstream inventory can help downstream units to lower safety stock.  Upstream nodes serving customers as well as multiple downstream facilities may need to prioritize between multiple orders.  Large but sporadic orders from downstream units can make it challenging for the upstream facilities to follow consistent inventory management policies, and thus it can lead to higher average inventory over time.  On the other hand, upstream facilities can choose a lower service level target as long as customer service levels from the downstream nodes are being met.  Every echelon can manage inventory with a custom policy instead of the standard reorder point or base stock approaches.  Modeling and optimization of such systems is therefore challenging as it requires a holistic approach that exploits such synergies and interactions between echelons \cite{chopra2007supply, evant2003}.

Several articles address inventory optimization for multi-echelon systems (see Axs{\"a}ter \cite{axsater1993continuous} for a review).  Based on the earliest model of Sherbrooke \cite{sherbrooke1968metric}, Ganeshan \cite{ganeshan1999managing} and Hopp et al.\ \cite{hopp1997easily,hopp1999easily} developed two-echelon models.  However, they assumed a reorder point policy, Poisson demand process for every echelon, and quantified expected backorders using infinite series.  Glasserman and Tayur \cite{glasserman1994stability,glasserman1995sensitivity} developed a multi-period model, but assumed a base-stock policy with 1 day review period and constant order placement.  Simpson \cite{simpson1958process}, Graves et al.\ \cite{graves2000optimizing, graves2003supply}, Inderfurth et al.\ \cite{inderfurth1991safety, inderfurth1998safety}, and You et al.\ \cite{you2010integrated} developed models for a base-stock policy using a \emph{guaranteed service time} approach in which each stage quotes a guaranteed service time.  For typical supply chains, where there are no quoted service time, this approach becomes very conservative with each stage planning for all previous stages' lead times and thus carry unnecessary excess inventory.  Lee and Billington \cite{lee1993material} and Ettl et al.\ \cite{ettl2000supply, ettl1999method} developed an algebraic \emph{stochastic service time} approach that characterizes random lead time delay induced downstream due to stockouts at upstream units.

Inventory optimization literature assumes a pre-defined probability distribution (such as Normal or Poisson), infinite series, or simplified custom algebraic expressions to quantify demand and lead time variability observed in supply chains.  Such assumptions can simplify or inaccurately represent the true nature of variability observed.  Moreover, most approaches make assumptions on system interactions, inventory policy, or the product flow in the system in order to be able to solve such problems with a mathematical programming framework.  We argue that this may not adequately capture the complex multi-echelon interactions described above and, as a results, can underestimate optimal inventory.  Therefore, in this work we follow a simulation-optimization framework to address this problem as it does not suffer from such assumptions (see \cite{agarwal2018validation}).

Major contributions of this work are:
\bit
	\item Instead of taking a parametric approach with a pre-decided probability distribution, the demand and the lead time variability are quantified by bootstrap sampling the historical data, thus preserving the true nature of the variability.  However, we do assume no time correlation in the historical data.  This can be easily relaxed by using a time-series model to forecast demand and lead time.
	\item We model the system using a flexible discrete-event simulation approach.  This allows us to capture sophisticated system interactions as well as any custom rules.
	\item The approach in this work is not limited to standard inventory policies such as the reorder point or base stock policies.  Any custom inventory policy can be incorporated.
	\item We compare three open source simulation optimization solvers: the derivative free methods from Scipy.Optimize \cite{jones2014scipy}, a Bayesian optimization algorithm Scikit-Optimize \cite{kumar2017scikit}, and a radial basis function based black-box optimizer RBFOpt \cite{costa2014rbfopt}.  There is no published work on this comparison.
\eit

The article is organized as follows.  The next section describes the problem addressed in this work and the example network.  Further in Section 3 we describe the inventory optimization algorithm with simulation and optimization modules.  We also illustrate a Python implementation of the algorithm in Section 4.  In Section 5 we present optimization results of the example network and a comparison of the open source solvers.  Section 6 concludes the article.

\section{Problem Description and Assumptions}

\begin{figure}[t]
	\centering
	\resizebox{13cm}{!}{\includegraphics[scale=1.0]{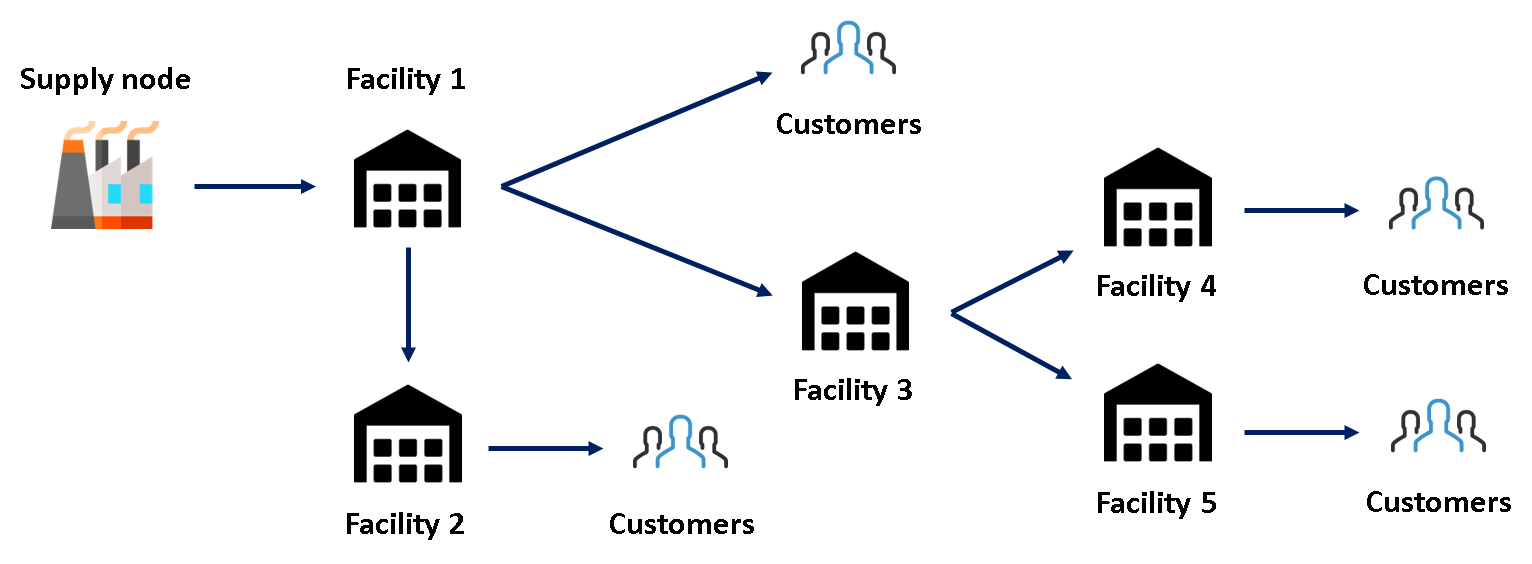}}
	\caption{A multi-echelon supply chain example}
	\label{fig:chain}
\end{figure}

Figure \ref{fig:chain} illustrates the multi-echelon supply chain network used in this work to demonstrate the simulation-optimization algorithm.  The network comprises both elements of a multi-echelon system where a stocking facility can either purely serve customer demand or purely replenishes another facility, or does a combination of both.  Here the outermost stocking facilities 2, 4, and 5 are customer facing units that directly serve customer demand.  Stocking facility 1 is a central facility that, besides directly fulfilling customer orders, also replenishes facility 2.  Facility 3 is another central facility that mainly replenishes 4 and 5.  It gets its replenishment from the first facility.  The first node that replenishes facility 1 is a supply node, such as a manufacturing plant or a vendor, for which we do not track inventory.

The purpose of inventory optimization is to minimize total system costs.  This can include order placement cost, inventory holding cost, and other miscellaneous charges such as facility costs, etc.  In this work we assume no order placement cost, and mainly focus on minimizing average on-hand inventory for all stocking locations while ensuring that a minimum customer service level target is met.  However, the framework presented is generic and any kinds of costs or parameters are straight-forward to incorporate.  We use $\beta$ service level (or fill rate) in this work \cite{agarwal2018validation}.  We choose to have a service level target for only customer serving facilities.  Thus, facility 3 has no explicit target; though it can centrally keep inventory for nodes 4 and 5 in order to facilitate meeting their service level targets.  The time gap between placing orders and receiving replenishment is characterized by an overall lead time.  Any unfulfilled customer order can either be back ordered or assumed to be lost sales; we study both cases.

Note that the facilities do not share their inventory levels or any other information with each other.  Each of them operate independently.  The optimization algorithm, being cognizant of the entire system, optimizes all echelons simultaneously.

\subsection{Key assumptions}
\bit
	\item Each stocking facility is replenished from only one upstream unit.
	\item All stocking units follow a combination of the reorder point and base stock policies.  When the inventory position falls below the reorder point, an order of amount of (base stock level - current inventory level) is placed by the facility.  We optimize both base stock and reorder point levels in this work.
	\item A replenishment order is never partially shipped.  If a stocking unit does not have enough on hand inventory, it waits until there is sufficient inventory to replenish the entire order volume.  Such a wait adds to the downstream unit's lead time.
	\item Serving customer order is prioritized over sending replenishment order.
	\item If multiple replenishment orders exist in the pipeline from multiple requesting facilities, they are fulfilled in the sequence they were received.
	\item Both customer demand and lead time exhibit variability.  Instead of assuming a pre-defined probability distribution, we follow a data-driven distribution to quantify variability.  We bootstrap samples from the historical data during simulation iterations.  This inherently assumes that there is no time correlation in the historical demand and lead time, as well as the future will be similar to the history.
\eit

\section{Inventory Optimization Algorithm}

We develop a simulation optimization framework to minimize inventory while ensuring a desired fill rate.  In this framework, the simulation module acts as a ``black-box'' for the optimization module.  The simulation process can accurately capture system dynamics.  Because of the black-box nature, the optimization model has visibility of the system dynamics only through the inputs to and outputs from the simulation model.  The optimization model invokes simulation module multiple times for function (and, if needed, derivative) evaluations.  

Note that, in general, the simulation-optimization frameworks do not guarantee optimality (or even feasibility) as many algorithms do not include a strong convergence criterion \cite{rios2012derivative}.  We ensure we obtain a feasible and a good quality solution by running the algorithm for several iterations.  This indeed results in a computationally slower framework.  However, we choose a simulation-optimization framework because:
\bit
	\item For practical inventory optimization problems, accurately capturing system interactions and dynamics is paramount
	\item A simulation model provides flexibility to model a variety of multi-echelon systems with custom business rules and inventory policies.  It can be easily tailored to the problem needs.
\eit

The next two sections describe the simulation and optimization algorithms in detail.

\subsection{Simulation Module}

The discrete-event simulation model comprises four different processes:
\ben
	\item \textbf{Place replenishment order:} Using this process the stocking locations place replenishment order to upstream facilities once their inventory levels reach below their reorder point .
	\item \textbf{Fulfill replenishment order:} This process is used by facilities to prepare and ship the replenishment order placed by their downstream facilities.  Once the order is prepared, it initiates the delivery process.
	\item \textbf{Deliver replenishment:} This process is invoked by the fulfill order process.  Each replenishment delivery is handled by this process.  It delivers the replenishment after waiting for downstream facility's lead time.  At delivery it increases downstream facility's on hand inventory.
	\item \textbf{Serve customer demand:} This process is used to deliver customer demand from each of the customer facing locations.
\een
The following sections describe each of these processes in detail.

\subsubsection{Place replenishment order}

\begin{procedure}
\DontPrintSemicolon

\SetKwFunction{Enqueue}{Enqueue}

\BlankLine
\KwIn{Facility $f$, Reorder point $R_f$, Base stock $B_f$, Upstream replenishing facility $S_f$, Simulation horizon $H$}
\KwResult{Inventory position $I_f$ and upstream order queue $Q_{S_f}$ get updated}
\BlankLine
\While{$t \leq H$}{
  \If{$I_{f,t} \leq R_{f}$}{
    $q \leftarrow B_{f} - O_{f,t}$\;
    \Enqueue{$Q_{S_{f}}, (q,f)$}\;
    $I_{f,t+1} \leftarrow I_{f,t} + q$\;
  }
  $t \leftarrow t + 1$\;
}
\caption{placeOrder($f$, $R_f$, $B_f$, $S_f$, $H$)}
\label{proc:placeOrder}
\end{procedure}

Procedure \ref{proc:placeOrder} illustrates the logic followed by each stocking location $f$ to order replenishment.  The procedure essentially tracks the inventory position $I_f$ of each facility.  As the inventory position falls below its reorder point $R_f$, an order is prepared that can bring the current on-hand inventory $O_f$ back to the base stock level $B_f$ (\emph{i.e.}, order volume = base stock - current on hand inventory).  The order is then inserted in the order queue of the upstream replenishing facility $Q_{S_f}$.  The inventory position is updated accordingly.

\subsubsection{Fulfill replenishment order}

Procedure \ref{proc:fulfillOrder} demonstrates the algorithm followed by stocking locations to ship the replenishment order placed by their downstream facilities.  As introduced in Procedure \ref{proc:placeOrder}, each facility maintains an order queue $Q_f$ which keeps track of the replenishment orders placed by its downstream facilities.  The queue follows a first-in-first-out mechanism.  Thus, downstream facilities' orders get processed in the order they were received.

Here, if the order queue $Q_f$ of facility $f$ is not empty, we retrieve the next order in the queue, which provides information on the order volume $o$ and its requesting facility $f'$.  It is ensured that the entire volume $o$ is shipped before the next order is considered from the queue.  Depending on the on-hand inventory volume $O_f$, we estimate if there is remaining volume $r$ that needs to be acquired in order to complete the order volume $o$.  If $r > 0$, the facility waits to ship the replenishment until there is enough on-hand inventory to fulfill $r$.  The on-hand inventory and inventory position are updated accordingly.  Once the order volume $o$ is available from the on-hand inventory this procedure calls the \ref{proc:shipOrder} procedure, described in the next section, that delivers the order to the requesting facility.

\begin{procedure}[t]
\DontPrintSemicolon

\SetKwFunction{Dequeue}{Dequeue}
\SetKwFunction{NotEmpty}{NotEmpty}
\SetKwFunction{Min}{Min}
\SetKwFunction{shipOrder}{shipOrder}

\BlankLine
\KwIn{Facility $f$, Simulation horizon $H$}
\KwResult{Inventory position $I_f$, on hand inventory $O_f$ and order queue $Q_f$ get updated.  Delivery of order $o$ for the downstream facility $f'$ is prepared}
\BlankLine
\While{$t \leq H$}{
  \eIf{\NotEmpty{$Q_f$}}{
    $(o,f') \leftarrow$ \Dequeue{$Q_f$}\;
    $s \leftarrow$ \Min{$o, O_{f,t}$}\;
    $O_{f,t} \leftarrow O_{f,t} - s$\;
    $I_{f,t} \leftarrow I_{f,t} - s$\;
    $r \leftarrow o - s$\;
    \If{$r > 0$}{
      \lWhile{$O_{f,t} < r$}{$t \leftarrow t + 1$}
      $O_{f,t} \leftarrow O_{f,t} - r$\;
      $I_{f,t} \leftarrow I_{f,t} - r$\;
    }
    \shipOrder{$f'$, $o$}\;
  }{
    $t \leftarrow t + 1$\;
  }
}
\caption{fulfillOrder($f$, $H$)}
\label{proc:fulfillOrder}
\end{procedure}

\subsubsection{Deliver replenishment}

Procedure \ref{proc:shipOrder} is symbolically identical to a conveyor belt.  Every replenishment delivery is conducted by this process.  When invoked by Procedure \ref{proc:fulfillOrder} it's primary function is to wait for the downstream facility $f$'s lead time before delivering the replenishment $o$.  At delivery it increases downstream facility's on hand inventory.

\begin{procedure}
\DontPrintSemicolon

\SetKwFunction{Bootstrap}{Bootstrap}

\BlankLine
\KwData{Historical lead time data $L_f$}
\KwIn{Facility $f$, Shipment volume $o$}
\KwResult{On hand inventory $O_f$ gets updated}
\BlankLine
$l \leftarrow$ \Bootstrap{$L_f$}\;
$t \leftarrow t + l$\;
$O_{f,t} \leftarrow O_{f,t} + o$\;
\caption{shipOrder($f$, $o$)}
\label{proc:shipOrder}
\end{procedure}

We note that the procedure utilizes historical lead time data for the requesting facility.  It bootstrap samples from the data to determine the lead time, instead of assuming any lead time probability distribution.  If the simulation time horizon is long enough, bootstrapping ensures that the nature of the lead time variability historically observed by $f$ in practice is also accurately captured and preserved in the simulations.

\subsubsection{Serve customer demand}

\begin{procedure}
\DontPrintSemicolon

\SetKwFunction{Bootstrap}{Bootstrap}
\SetKwFunction{Min}{Min}
\SetKwFunction{Max}{Max}

\BlankLine
\KwData{Historical customer demand $D_f$}
\KwIn{Facility $f$, Simulation horizon $H$, $choice$ of lost sales or back order}
\KwResult{Inventory position $I_f$, on hand inventory $O_f$, and total demand $M_f$ get updated. If $choice$ is lost sales, total shipped volume $P_f$ gets updated. If $choice$ is back order, back order volume $B_f$ and total back ordered volume $T_f$ get updated.}
\BlankLine
\While{$t \leq H$}{
  $d \leftarrow$ \Bootstrap{$D_f$}\;
  $M_{f,t+1} \leftarrow M_{f,t} + d$\;
  \eIf{$choice$=Lost sales}{
    $s \leftarrow$ \Min{$d, O_{f,t}$}\;
    $P_{f,t+1} \leftarrow P_{f,t} + s$\;
  }{
    $s \leftarrow$ \Min{$d + B_{f,t}, O_{f,t}$}\;
    $b \leftarrow d - s$\;
    $B_{f,t+1} \leftarrow B_{f,t} + b$\;
    $T_{f,t+1} \leftarrow T_{f,t}$ + \Max{$0, b$}\;
  }
  $O_{f,t+1} \leftarrow O_{f,t} - s$\;
  $I_{f,t+1} \leftarrow I_{f,t} - s$\;
  $t \leftarrow t + 1$\;
}
\caption{serveCustomer($f$, $H$, $choice$)}
\label{proc:serveCustomer}
\end{procedure}

Finally, we describe the procedure that fulfills customer demand from the customer facing stocking locations.  Procedure \ref{proc:serveCustomer} incorporates both cases when demand is either back ordered or considered as lost sales.  The algorithm follows either of the two cases depending on a pre-specified variable $choice$.

Here, as in the Procedure \ref{proc:shipOrder}, we capture demand variability by bootstrapping samples from the historical customer demand data observed by the serving facility.  Thus, we ensure accurate demand variability is captured and preserved in the simulations.  Note that the demand data is at the granularity of the simulation time period.  For instance, if simulating daily inventory behavior, the historical data should be aggregated daily demand.  Also, we assume there is no time correlation in the historical demand in order to randomly sample from the data.

Procedure \ref{proc:serveCustomer} tracks total demand observed $M_f$ in order to estimate $\beta$ service level (described in the next section).  First we consider the lost sales case. Here if customer demand $d$ can be met with the on-hand inventory $O_f$, it is shipped, otherwise the volume available on-hand is delivered.  We record total shipment delivered $P_f$.  $P_f$ together with $M_f$ are used to estimate $\beta$ service level.

In the case when the unmet demand can be back ordered, we first determine if there is enough on-hand inventory to meet current period demand and previous back orders.  If yes, the entire volume is shipped, otherwise the volume available on-hand is delivered.  We record the difference between demand $d$ and shipped volume $s$ in $b$.  Note that if we ship previous back orders along with $d$, $b$ is negative.  $b$ is used to update the back orders tracker $B_f$.  Thus, a negative $b$ reduces the volume of back orders.  If $b$ is positive, it gets recorded as late sales $T_f$.  $T_f$ along with $M_f$ are used to estimate service level for the case of back orders (described in the next section).

Depending on the shipment volume we accordingly update the on-hand inventory $O_f$ and inventory position $I_f$.

\subsubsection{Simulation initialization}

\begin{procedure}[t]
\DontPrintSemicolon

\SetKwFunction{placeOrder}{placeOrder}
\SetKwFunction{fulfillOrder}{fulfillOrder}
\SetKwFunction{serveCustomer}{serveCustomer}

\BlankLine
\KwIn{List of facilities $F$, Reorder point $R_f$, Base stock level $B_f$, and upstream replenishing facility $S_f$ for all facilities $f \in F$, simulation horizon $H$, and $choice$ of lost sales or back order}
\KwResult{Estimate average on-hand inventory $A_f$ and $\beta$ service level $\beta_f$ for $f \in F$ }
\BlankLine
\ForEach{$f \in F$}{
  \placeOrder{$f$, $R_f$, $B_f$, $S_f$, $H$}\;
  \fulfillOrder{$f$, $H$}\;
  \serveCustomer{$f$, $H$, $choice$}\;
}
\ForEach{$f \in F$}{
  $A_f \leftarrow 0$\;
  \lWhile{$t \leq H$}{$A_f \leftarrow A_f + O_{f,t}$}
  $A_f \leftarrow A_f/H$\;
  \eIf{$choice$=Lost sales}{
    $\beta_f \leftarrow P_{f,H}/M_{f,H}$\;
  }{
    $\beta_f \leftarrow 1 - T_{f,H}/M_{f,H}$\;
  }
}
\caption{simNetwork($F$, $R$, $B$, $S$, $choice$, $H$)}
\label{proc:simulation}
\end{procedure}

Procedure \ref{proc:simulation} is the main initialization function of the simulation module.  It connects all procedures described above.  It initiates simulation and computes average on-hand inventory and $\beta$ service level or fill-rate for each facility.  The optimization module invokes this procedure for function and derivative evaluation.

When invoked, Procedure \ref{proc:simulation} initiates procedures \ref{proc:placeOrder}, \ref{proc:fulfillOrder}, and \ref{proc:serveCustomer} for each facility, which run in parallel for the entire simulation horizon $H$.  It requires list of facilities $F$, reorder point $R_f$, base stock level $B_f$, the supply chain network and upstream replenishing facility $S_f$ for all facilities $f \in F$, simulation horizon $H$, and $choice$ of lost sales or back order as inputs, which are further passed to other procedures.  Note that during the optimization process, some of these inputs can be decision variables and are provided by the optimization module.

Once the simulation is complete, the procedure computes average on-hand inventory $A_f$ by taking a sample average across $H$.  Also, as detailed in the previous section, $\beta$ service level is obtained using $P_f$, $T_f$, and $M_f$ depending on the $choice$ of considering unfulfilled demand as lost sales or back ordered.  $A_f$ and $\beta_f$ are utilized by the optimization module in its objective to determine optimal decisions.

\subsection{Optimization Module}

Procedure \ref{proc:optimization} illustrates the module that is invoked by a black-box optimization solver for function and derivative evaluation.  It computes the objective function for the current value of decision variables provided by the black-box solver.  

The objective function has two components that the solver seeks to minimize.  The first component is the total of average on-hand inventory of all facilities, while the second is the penalty for not meeting the service level targets summed across all facilities.  Both components are computed for multiple simulation replications and averaged across all replications to create a composite objective function for the optimization solver.   

In each replication, \ref{proc:optimization} invokes the \ref{proc:simulation} procedure to obtain average on-hand inventory $A_f$ and fill rate $\beta_f$ for each facility $f$.  The current value of the inputs ($F$, $R$, $B$, $S$, $choice$, $H$) obtained from the black-box optimization solver are provided to \ref{proc:simulation}.

The average on-hand inventory is aggregated across all facilities in $AA$. The fill rate $\beta_f$ is compared against the target service level $\beta_{f}^{T}$.  Any unfulfilled service level $\beta_{f}^{T} - \beta_f$ is also accumulated in $A\beta$.  This is repeated for user-specified replications $N$. Finally, $AA$ and $A\beta$ are averaged across $N$, and $A\beta$ is penalized with a user-specified violation penalty $\rho$ to obtain the objective function $Z$.

Note that, as mentioned before, we do not include order placement or any other costs in this framework.  However, if needed, they can be easily added to the algorithm as additional terms in the objective function and computed in a similar manner as demonstrated for average on-hand inventory and fill rate.

\begin{procedure}
\DontPrintSemicolon

\SetKwFunction{simNetwork}{simNetwork}
\SetKwFunction{max}{max}

\BlankLine
\KwData{Number of simulation replications $N$, target service level $\beta^{T}_{f}$ for each facility $f$, service level violation penalty $\rho$}
\KwIn{List of facilities $F$, Reorder point $R_f$, Base stock level $B_f$, and upstream replenishing facility $S_f$ for all facilities $f \in F$, simulation horizon $H$, and $choice$ of lost sales or back order}
\KwResult{Determine objective function $Z$ for the black-box optimization algorithm based on the average of the average on-hand inventory $AA$ and $\beta$ service level $A\beta$ across all replications $N$}
\BlankLine
$AA \leftarrow 0$\;
$A\beta \leftarrow 0$\;
\ForEach{$n \in N$}{
  $(A_f,\beta_f) \leftarrow $ \simNetwork{$F$, $R$, $B$, $S$, $choice$, $H$}\;
  \ForEach{$f \in F$}{
    $AA \leftarrow AA + A_f$\;
    $A\beta \leftarrow A\beta + $ \max{$0,\beta^{T}_{f} - \beta_f$} \;
  }
}
$Z \leftarrow AA/N + \rho A\beta/N$\;
\caption{optNetwork($F$, $R$, $B$, $S$, $choice$, $H$)}
\label{proc:optimization}
\end{procedure}

\section{Algorithm Implementation}

\subsection{Optimization Solvers}

As mentioned above, the optimization module of the algorithm is invoked by a simulation-optimization solver, which determines the optimal decision variables.  In this work we utilize and compare three open source libraries: the derivative free methods from Scipy.Optimize \cite{jones2014scipy}, a Bayesian optimization algorithm Scikit-Optimize \cite{kumar2017scikit}, and a radial basis function based black-box optimizer RBFOpt \cite{costa2014rbfopt}.  Below we provide a brief note on the methodologies used by these solvers.

\paragraph*{SciPy.Optimize}  This is a set of optimization algorithms provided by the scipy Python package.  It provides various algorithms for both constrained and unconstrained optimization.  The general purpose algorithms for unconstrained optimization include Nelder-Mead, Powell, Conjugate Gradient (CG), BFGS, and Newton-CG.  Constrained optimization algorithms include L-BFGS-B, truncated Newton, COBYLA, SLSQP, and various trust-region methods.  Quasi-Newton strategies SR1 and BFGS are available to approximate Hessian update.  Several algorithms require gradient or Jacobian information to be provided.  In this work we use the derivative free algorithms.

\paragraph*{Scikit-Optimize or skopt} This is a non gradient-based optimization algorithm that minimizes expensive and noisy black-box functions.  It is built on NumPy, SciPy, and Scikit-Learn Python packages.  It essentially utilizes a Bayesian optimization approach.  Here the unknown objective is considered as a random function for which a prior distribution is defined using a Gaussian process.  Function evaluations are treated as data and used to update the prior to form the posterior distribution for the objective function.  This posterior is used to maximize a simple utility function.  Scikit-Optimize also provides capabilities to, instead of a Gaussian process, use Random Forest and XGBoost algorithms to approximate the noisy functions.

\paragraph*{RBFOpt} This library implements the Radial Basis Function method originally proposed by Gutmann \cite{gutmann2001radial} with a few extensions.  Here they build and iteratively refine a surrogate model of the unknown objective function.  First, for an initial set of chosen sample points that satisfy the bounds and integrality constraints, a surrogate interpolation function is built using a combination of radial basis functions and polynomials.  The accuracy of the surrogate model is assessed and model selection is performed automatically in the algorithm using a cross validation scheme.  Next step is to choose a trade-off between \emph{exploration} and \emph{exploitation} based on Gutmann's idea of ``bumpiness''.  \emph{Exploration} implies trying to improve the surrogate model in unknown parts of the domain, whereas \emph{exploitation} implies trying to find the best objective function value based on the current surrogate model.  The next point in the search space, to which the algorithm moves, is determined based on this trade-off.  Gutmann's algorithm provides strategies to evaluate this trade-off based on least ``bumpy'' points, and solves a local search optimization problem to determine the next iterate.

\subsection{Implementation Code}

An open-source Python-based implementation of our inventory optimization algorithm is available here \cite{github}.  We implement the simulation framework using the open-source SimPy platform \cite{matloff2008introduction}.  Both back order and lost sales simulations are implemented separately.  A user needs to pre-specify which simulation option is used by the optimization routines.  For optimization, we provide separate implementation for all three black-box optimization libraries.  We compare the results obtained with the three solvers in the next section.

\section{Results and Solver Comparison}

We demonstrate the inventory optimization algorithm and compare the open source optimization libraries on the example supply chain network shown in Figure \ref{fig:chain}.  As described before, all facilities follow a combination of the reorder point and base stock policies.  When the inventory position falls below reorder point, an order of amount = (base stock level - current inventory level) is placed by the facility.  We optimize both base stock and reorder point levels.

\begin{table}[b]
\begin{center}
\caption{Parameters for the supply chain in Figure \ref{fig:chain}}
\begin{tabular}{|c|c|c|c|c|}
\hline
\T Facility & Base & Service level & Base stock & Reorder point \\
\B & lead time & target & initial guess & initial guess \\
\hline
\T Facility 1 & 3 days & 95\% & 3000 & 1000 \\
Facility 2 & 4 days & 95\% & 600 & 250 \\
Facility 3 & 4 days & 0\% & 900 & 200 \\
Facility 4 & 2 days & 95\% & 300 & 150 \\
\B Facility 5 & 2 days & 95\% & 600 & 200 \\
\hline
\end{tabular}
\label{tbl:params}
\end{center}
\end{table}

Table \ref{tbl:params} shows the base lead time, service level target, and initial guesses for base stock level and reorder point assumed for all facilities in the supply chain.  The initial guess is feasible, \emph{i.e.}, the facilities satisfy their respective service level target with this initial guess.

The base lead time is the minimum replenishment lead time from the corresponding supplier facility.  The true lead time experienced is the base lead time + a random variable component.  This random variable component is obtained from the historical data.  During the simulation-optimization, we bootstrap from the historical data to get this variable component.  Similarly, the customer demand for each facility is also bootstrapped from the historical data during simulations.  The historical data assumed and used for this example is not correlated with time.

Note that in Table \ref{tbl:params}, the service level target for Facility 3 is zero.  This is not impractical.  Facility 3's role predominantly is to support replenishment for Facilities 4 and 5.  It need not have a service level target as long as it is able to help Facilities 4 and 5 meet their service level targets, which are 95\%, respectively.  While centralized inventory at Facility 3 can pool the demand fluctuation risk at 4 and 5, decentralized individual inventory at both Facilities 4 and 5 with a lower lead time can lead to a more responsive supply chain.  The service level target for Facility 1 applies only to the customer demand; replenishment to Facility 2 is not accounted in calculating its service level.

For the objective function we use a penalty of $10^6$ for not meeting the service level target.  We perform 20 simulation replications and take an average across replications for the average on-hand inventory and service level penalty for the objective function.  We assume the initial inventory to be 90\% of the base stock level for all facilities.  Also, each simulation is run for a 360 days time period.  Both simulation and optimization are performed on a Mac OS with a 1.2 GHz Intel Core m3 processor and 8 GB 1867 MHz RAM.

Table \ref{tbl:algoparams} shows the black box algorithm used from each Python package and their parameter settings.  Below we describe our implementation strategy of each algorithm:
\bit
	\item{\textbf{SciPy.Optimize:}} Here we run the Nelder Mead algorithm for 100 cycles.  In each cycle the algorithm runs for 50 iterations.  After 50 iterations another cycle starts with the current best solution as the initial guess for this new cycle.  The cycle then runs for another set of 50 iterations.  This sequence continues until we either exceed the maximum time limit of 1 day or complete all 100 cycles.
	\item{\textbf{Scikit-Optimize:}} This algorithm relies heavily on the starting random state.  Instead of the \emph{exploitation} approach, where we run the algorithm longer for several iterations for one particular starting random state, in our tests the algorithm performed significantly better with the \emph{exploration} approach, where the algorithm is run for multiple various starting random states with fewer iterations per run.  We run the algorithm for 1000 cycles, each cycle explores a different start state with 20 iterations.  We varied the random start state from 0 to 1000.  A Kappa of 50 was chosen to further enhance the \emph{exploration} capability of the algorithm.  The procedure was run until we either exceeded the maximum time limit of 1 day or completed all 1000 cycles
	\item{\textbf{RBFOpt:}} Because of a thorough and robust internal implementation of the algorithm that explores various strategies to search the space, we simply ran the algorithm as is for 1000 iterations.  We chose a random seed of 707 in order to be able to reproduce the results.
\eit

\begin{table}
\begin{center}
\caption{Settings and parameters for the optimization solvers}
\begin{tabular}{|c|c|l|}
\hline
\T Python Package \B & Algorithm & Settings \\
\hline
\T SciPy.Optimize & Nelder Mead & Number of cycles = 100 \\
& & Iterations per cycle = 50 \\
\B & & Max time limit = 1440 minutes \\
\hline
\T Scikit-Optimize & gp\_minimize & Number of cycles = 1000 \\
& & Iterations per cycle = 20 \\
& & Number of random starts = 10 \\
& & Random state $\in$ [0,1000] \\
& & Kappa = 50 \\
\B & & Max time limit = 1440 minutes \\
\hline
\T RBFOpt & RbfoptUserBlackBox & Max iterations = 1000 \\
& RbfoptAlgorithm & Global search method = \emph{solver} \\
\B & & Random seed = 707 \\
\hline
\end{tabular}
\label{tbl:algoparams}
\end{center}
\end{table}

\begin{table}[t]
\begin{center}
\caption{Optimal results and solver comparison when the unmet demand is back ordered}
\begin{tabular}{|l|c|c|c|}
\hline
\T \B & SciPy.Optimize & Scikit-Optimize & RBFOpt \\
\hline
\T Final objective & & & \\
\hspace{3mm} Optimal objective & 2516 & 1517 & 951 \\
\B \hspace{3mm}\% reduction from the initial guess & 7\% & 44\% & 65\% \\
\hline
\T Optimal Base stock & & & \\
\hspace{3mm}Facility 1 & 2736 & 2034 & 835 \\
\hspace{3mm}Facility 2 & 635 & 383 & 254 \\
\hspace{3mm}Facility 3 & 950 & 527 & 418 \\
\hspace{3mm}Facility 4 & 307 & 190 & 181 \\
\B \hspace{3mm}Facility 5 & 598 & 337 & 246 \\
\hline
\T Optimal ROP & & & \\
\hspace{3mm}Facility 1 & 804 & 709 & 835 \\
\hspace{3mm}Facility 2 & 257 & 187 & 196 \\
\hspace{3mm}Facility 3 & 213 & 181 & 188 \\
\hspace{3mm}Facility 4 & 149 & 91 & 84 \\
\B \hspace{3mm}Facility 5 & 199 & 188 & 200 \\
\hline
\T Other diagnostics & & & \\
\hspace{3mm}Total iterations & 5,000 & 20,000 & 1,000 \\
\B \hspace{3mm}CPU time (minutes) & 248 & 399 & 94 \\
\hline
\end{tabular}
\label{tbl:solnbackorder}
\end{center}
\end{table}

\begin{table}[t]
\begin{center}
\caption{Optimal results and solver comparison when the unmet demand is considered as lost sales}
\begin{tabular}{|l|c|c|c|}
\hline
\T \B & SciPy.Optimize & Scikit-Optimize & RBFOpt \\
\hline
\T Final objective & & & \\
\hspace{3mm} Optimal objective & 2446 & 1277 & 1146 \\
\B \hspace{3mm}\% reduction from the initial guess & 12\% & 54\% & 59\% \\
\hline
\T Optimal Base stock & & & \\
\hspace{3mm}Facility 1 & 2516 & 1854 & 1550 \\
\hspace{3mm}Facility 2 & 643 & 252 & 281 \\
\hspace{3mm}Facility 3 & 937 & 424 & 451 \\
\hspace{3mm}Facility 4 & 308 & 124 & 99 \\
\B \hspace{3mm}Facility 5 & 625 & 190 & 259 \\
\hline
\T Optimal ROP & & & \\
\hspace{3mm}Facility 1 & 729 & 732 & 722 \\
\hspace{3mm}Facility 2 & 276 & 195 & 214 \\
\hspace{3mm}Facility 3 & 198 & 189 & 135 \\
\hspace{3mm}Facility 4 & 159 & 106 & 99 \\
\B \hspace{3mm}Facility 5 & 220 & 182 & 200 \\
\hline
\T Other diagnostics & & & \\
\hspace{3mm}Total iterations & 5,000 & 20,000 & 1,000 \\
\B \hspace{3mm}CPU time (minutes) & 237 & 387 & 117 \\
\hline
\end{tabular}
\label{tbl:solnlostsales}
\end{center}
\end{table}

Table \ref{tbl:solnbackorder} illustrates the optimization results when we consider the unmet demand to be back ordered, while Table \ref{tbl:solnlostsales} shows the optimization results when the unmet demand is considered as lost sales.  For all three solvers, we compare the optimal base stock and reorder points, total number of iterations and computational time, final value of the objective function, and how much reduction in the objective was achieved from the corresponding value for the initial guess.  Note that the CPU time is based on 20 simulation replications per optimization iteration.  For more replications the CPU time will increase.  Also, the optimal base stock and reorder points pertain to the example data used in the case study.  We include them in the table in order to perform a comparative analysis between the solvers.  We observe the following from the results:
\bit
	\item Among the three solvers, RBFOpt produced the best results for both back order and lost sales cases.  It was able to achieve the best reduction in the objective function and the lowest inventory levels.  Moreover, it was able to achieve results in fewest iterations and within the best CPU time of less than two hours among the three solvers.
	\item SciPy.Optimize with the Nelder Mead algorithm doesn't seem to perform well.  Even when ran for 5000 iterations for more than 4 hours, it could only achieve 7\% reduction in the objective function in the back order case, and 12\% reduction in the lost sales case, respectively.  The final optimal ROP and base stock are quite high and not very different from the initial guess.
	\item While Scikit-Optimize produced considerably better results than SciPy.Optimize, it took four times more iterations and 40\% more CPU time.  It is computationally most expensive compared to the other two because of its overall methodology.  Because the algorithm relies more on the random starting state, while it can produce better results, the exploration approach with multiple random start states is computationally more intensive.
	\item The performance of Scikit-Optimize is dependent on what starting states are searched.  Tables \ref{tbl:solnbackorder} and \ref{tbl:solnlostsales} illustrate the best results obtained among various tests with different random start states for the case study in this work.  For a different problem, the statistics could be different.  For instance, while RBFOpt appears to be a clear winner when compared to Scikit-Optimize for the back order case, the final optimal objective for both are quite close for the lost sales case.  Therefore, it is not definitive that Scikit-Optimize is better or worse than RBFOpt.  Depending on the problem, search strategy, and random start states, good solution can be obtained from either of the two algorithms.  However, in all our tests, RBFOpt's results were always superior compared to the others.
	\item While we observe either comparable or higher optimal ROP from RBFOpt when compared to that from Scikit-Optimize, the base stock values are lower, especially when observed for Facility 1, thus leading to an overall lower average inventory.  The optimal ROP values are almost comparable across all three solvers.  Both Scikit-Optimize and RBFOpt are able to achieve lower inventory by attaining lower base stock values when compared to SciPy.Optimize.
	\item In general, when systems consider unmet demand as lost sales instead of back ordering them, they tend to carry lower inventory on average with lower reorder points.  We observe this trend in our case studies with both SciPy.Optimize and Scikit-Optimize.  However, RBFOpt in its results show the opposite behavior.  Also, while SciPy.Optimize and Scikit-Optimize were faster in the lost sales case compared to back order, RBFOpt took longer to complete iterations.  This indicates potential for getting a more improved solution with RBFOpt for the lost sales case if it is run longer for more iterations.
	\item For both back order and lost sales scenarios, we observe that the algorithms trend towards keeping a lower ROP at the customer facing facilities while compensating that with a higher ROP at intermediate warehouses.  This behavior is more apparent in the lost sales case.  We infer that optimal systems demonstrate a preference for a centralized inventory planning and positioning scheme that help with risk pooling.
	\item We note that the optimal ROP for Facility 5 from RBFOpt is same as the initial guess.  This indicates potential for a further improved solution from RBFOpt if it is run longer for more iterations
	\item One key observation is that, for Facility 1 in the back order case, and for Facilities 2 and 4 in the lost sales case, RBFOpt's solution kept base stock same as the ROP.  This is the influence of assuming no order placement cost for our case study as the system now prefers to place frequent orders with no additional penalty.  Thus, it is important to make appropriate assumptions to model a system as it significantly impacts the optimal solution and the behavior. 
\eit

\section{Conclusions}

We approach the problem of optimizing inventory in multi-echelon supply chain systems by developing a simulation-optimization framework.  The simulation capability, by design, can model sophisticated system interactions as well as any custom rules.  Moreover, a black-box optimization algorithm wrapped around the simulation allows determining optimal system planning decisions while ensuring system dynamics are accurately captured.  In order to preserve the true nature of the variability experienced by the system, we quantify demand and lead time variability by bootstrapping historical data.

In order to demonstrate practical applicability of the framework, we test it on a 3-echelon example network.  Average inventory is minimized while ensuring that the desired average $\beta$ service level at stocking locations is achieved.  We assume no order placement cost.  Because there is no single proven simulation optimization platform, we compare three open source solver packages: Scipy.Optimize, Scikit-Optimize, and RBFOpt.

The key conclusions from the work are as follows:
\bit
	\item While we observe substantially lower inventory levels and computationally superior results from RBFOpt, depending on the problem, search strategy, and the random start states, close enough good solution can be obtained from both RBFOpt and Scikit-Optimize.
	\item The optimization results demonstrate a preference for a centralized inventory planning scheme that help with risk pooling.
	\item Because we assume no order placement cost, as expected in the optimal solution, the facilities tend to order more frequently in order to lower their inventory as there is no additional penalty.  Thus, it is important to make appropriate assumptions to model a system as it significantly impacts the optimal solution and the modeled behavior.
\eit

\bibliographystyle{amsplain}
\bibliography{inventory_opt_refs}

\end{document}